# The beginnings of the Soviet encyclopedia.

## *Utopia and misery of mathematics in the political turmoils of the 1920s.*

Laurent MAZLIAK[1]

**Abstract :** *In this paper we focus on the beginning of publication of the Large Soviet Encyclopedia, launched in 1925. We present the context of this launching and explain why it was tightly connected to the period of the New Economical Policy. In a last section, we examine four articles included in the first volumes of the encyclopedia and relative to randomness and probability, in order to illustrate some debates of the scientific scene in USSR during the 1920s.*

**Key-Words and Phrases :** *National encyclopedia. USSR. New Economical Policy. History of probability theory. Foundations of Probability. Otto Schmidt. Arthur Bowley.*

**Introduction**

The Large Soviet encyclopedia (Большая Советская Энциклопедия ; LSE in the sequel) was a gigantic enterprise to the glory of "Marxist science" and of the Soviet regime. There were three editions : the first one was launched in 1926, the second one in 1949, the third one in 1977. The statement of these years alone allows to understand the enormous differences between the three editions. Moreover, one observes that if the second and the third were launched in a relatively short period of time, nine years for both, the first edition needed more than twenty years to be completed.

Roughly speaking, the second edition of the LSE was characterized by the years of Stalinist glaciation of the 1950s, when, after World War 2, Soviet Union and its satellite countries were more or less isolated behind the iron curtain. As for the third one, it represented the last attempt for the declining regime to present a general picture of the Soviet conception of the world and it implied a deep tidying up of the most salient aspects inherited from the Stalinist period. Thus, clearly, the first edition represents the richest of the three editions as historical source for a better understanding of how Soviet thinking was constructed after 1917. There are several reasons why this edition provides such a capital wealth. First of all, it is so precisely because it was the first edition : it imposed several forms to the publication which would be continued in the following editions. Among these forms, the most obvious, which bore innumerous consequences, was the choice of the alphabetical ordering for the entries, so that the publication can be seen as an encyclopedic dictionnary as well as an encyclopedia. Some other choices were kept in future editions : a rather small dimension for the volumes, a two-column display of the pages, two sizes of fonts with large or small letters. Also the presence of numerous pictures and drawings. Another reason for this edition to be precious is that it included among its collaborators a huge number of first-rate personalities of the Soviet academic scene of the time. But the most important of all the reasons for this first edition to be so important is in fact that the publication of its 55 volumes lasted more that 20 years, between 1926 and 1947, and therefore it witnessed the enormous changes met by Soviet

---

[1] LPMA, University Pierre and Marie Curie, Paris, France. laurent.mazliak@upmc.fr

Union during this period. Even if we limit ourselves to the first five years, as we shall do in the present paper, it is worth recording that the 1920s in USSR began in the violence of war communism, followed since 1922 by the very particular time of the New Economical Policy (НЭП in Russian - NEP in the sequel). The decade was abruptly closed after 1928 by the end of the NEP and, in the academic world, by the fight against bourgeois specialists, a herald for the nightmare of the 1930s. The volumes of the LSE published during, say, the five years 1926-1930, joigning the immediate aftermath of Lenin's death (1924) to the consolidation of the Stalinist dictatorship were written on the background of the sinuosities of Soviet orthodoxy during this period and are precious to document the history of these difficult times, especially the intellectual history, marked by the still vivid ambition that Bolshevik Russia would be the spearhead of a world revolution.

From the very origin, an ideological basis was searched for the project of publishing an encyclopedia (we shall come back later at length in the paper on the context of this project), in order to justify its necessity and its coherence with the educational propaganda of the regime. The promoters of the project were interested in proving that in the past, the importance of this kind of enterprise had already been mentioned in relation to the proletarian revolutionnary movement. For instance, the following quotation from a text written by Jaurès in 1901 was used for that purpose[2]

> *... In my eyes the hour comes closer when the socialist and revolutionary proletariat must acquire an organized doctrine of the universe and of life. What the Encyclopedia has been for the revolutionary bourgeoisie, a new encyclopedia, infinitely bolder and wider, will have to be for the proletariat. We shall have to resume the movement of human thought from Kant to Renan, through Hegel, Comte and Marx. We shall have to resume the movement of science from Laplace to Maxwell, through Darwin, to offer the key findings and the main trends to the proletariat who wants to live its life to the full, and to project a bright light on the universe where enlightments of individual thinking will mix with the fiery radiance of social life [...] There must be a general philosophy, both revolutionary and evolutionary, which is gradually communicated to the conscious proletarian elite, and by degrees to the whole proletariat.* (Jaurès,1901)

It is seen that the use of the previous quotation was possible only by a slight shift of its original context. It is indeed not quite clear that Jaurès employed the word "encyclopedia" with a concrete publication enterprise in mind. When the first volume was published in 1926, the editorial board shelled out a preface to expose the general program of the LSE, for which the scientific method was presented as the very principle of the book because it perfectly suited the political aims of the new regime. As the communist economist Maria N.Smit-Falker (to which we shall come back at length in the third part of the present paper) expressed

> *The more the socialism in our country will progress, the stronger will be the influence of the scientific thinking on life, and the greater will be the role of the scientific and social organisations for the resolution of practical*

---
[2] Quoted in (Петров, 1960) p. 132.

*questions.* (Smit-Falker, 1927; p.15)

A major theme for the aforementioned preface from the editorial board was the thorough transformation of the political situation in Russia which had implied the emergence of a new kind of readership.

> *The Revolution created a new reader, with new questions, with the persistent desire to get an orientation in all the variety of the contemporary world, to systematize his knowledge, to strengthen his conception of the revolutionary and materialistic world, to get acquainted with the last advances of science. Our time is a period of transition of capitalism towards socialism, when in a fundamental way are transformed the material bases as well as the social relationships and the ideology.*

This argument, in turn, was used to justify the need for a new kind of publication - and in this case, a new kind of encyclopedia. The LSE had to be of a different nature from the previous encyclopedias. The board was eager to emphasize the differences between former enterprises (especially those from the Russian czarist-period) and the new publication, due in particular to a different ideological approach.

> *In previous dictionaries, coexisted different conceptions of the world - even contradictory. On the contrary, for the Soviet encyclopedia, a precise vision of the world is absolutely necessary, namely a strictly materialistic one. Our worldview is dialectical materialism. Humanities, both to understand the past and modern times, have already been extensively transformed on the basis of a continued application of the dialectical method of Marx and Lenin; in natural and exact sciences, the board, while trying to highlight the standpoint of dialectical materialism, will take account of the fact that there is not in all areas a sufficient number of perfectly Marxists studies. In these sciences is hardly built the foundation required for the implementation of the dialectical method. The encyclopedia promotes strictly factual side of natural sciences, released from their idealistic premises.*

Naturally, it was also necessary to make some concessions to the spirit of the time and to assert that the volumes of the LSE would be accessible to the alleged new master of the Soviet society, the factory workers.
The continuation of the preface answers to that concern with a noticeable shift of the "center of gravity" of the topics towards practical application, social and political construction and a more or less accepted side-lining of abstraction.

> *In the previous dictionaries one felt that they were written for scholars with an interest primarily in literature and history - on the contrary the LSE has translated to the social sciences the center of gravity: in economy, contemporary politics and Soviet practice. To exact and natural sciences is attributed a large place, but not for the dry description of the different kinds of plants or various abtract questions. The natural sciences for the LSE are the foundation of the work for the domination over the forces of nature and for their use for human needs. Therefore a much more important place than previously is allocated to agriculture, industry and*

> *technology. At the center of our attention : the Soviet Union, the construction of our society and governance and the international revolutionary movement.*

One finds in this tirade some elements directly inherited from the rationalist movements of the turn of the 20th century. It is striking how the mention of the scientific method in the previous quotations is in line with how the mathematician Emile Borel, in 1906, presented the scope of his newly founded journal, the *Revue du Mois*. One reads in the foreword of the first issue of the journal

> *The number and the importance of problems that can be treated by adopting scientific methods grows every day. It seemed possible to us to imagine a journal which focused on these methods, not as a specialist publication but rather by aiming at the general development of ideas, and the exposition and critical appraisal of the advances in Knowledge and the resultant spread of ideas.*
> *The Revue du Mois attempts to be this journal. It claims, above all, to be a journal containing free discussion, allowing the free unhampered expression of opinions based on science.* (Borel, 1906)[3]

The concept of scientific method was naturally not exactly the same in Borel's and in the editors of the LSE's mind. For the latters, the scientific method was to be tested through the sieve of a "Marx-Engelsisation" with a Leninist touch : dialectical materialism. Observe nevertheless that the editorial board of the LSE carefully provided an alibi in the preface by admitting that the situation was not uniform over all the topics because for some of them, especially the natural and exact sciences, Marxist science was not enough advanced. It clearly followed that in these domains it was "classical bourgeois" science that would generally be exposed in the LSE. The board obviously had the wisdom not to write down explicitly this consequence.

The aim of the present paper is to give information about the first years of the Large Soviet encyclopedia and particularly about some aspects of mathematics in it, with a special focus on probability theory as an illustration of the ideological background of the publication. The paper is organized as follows. In the first part, we comment on the Soviet conception of the scientific method, and how the intelligentsia faced the new regime's demands. In a second part, we present the origins of the LSE project, its responsibles and its connection with the Soviet politics of the 1920s. In the third and last part, we focus on the question of the mathematics of randomness in USSR and illustrate this question by four entries belonging to the first volumes of the LSE.

I - **Soviet ideology and scientific method**

1- Some comments on dialectical materialism in USSR

In his classical book (Graham, 1987) (especially chapter 2, p.25-67), from which we take the elements of the present subsection, Loren Graham describes in detail how the concept of dialectical materialism became central to Soviet thinking. This concept was used for the first time under that name by Plekhanov, the "father of Russian Marxism" in 1891 in a comment

---
[3] For details about Borel's foundation of the Revue du Mois, the interested reader may consult (Durand and Mazliak, 2011).

on Hegel. However, it was Engels' writings which were at the source of the scientific practice in Soviet dogma. In the preface of his anti-Dühring, Engels writes that "a knowledge of mathematics and natural sciences is essential to a conception of nature which is dialectical and at the same time materialist" (Engels, 1959; p.16)

Engels's idea was that the aim of such a knowledge allows to concentrate on general laws describing the process acting in the material world. The notion of law (*Gesetz*) for Engels was rather vague, and especially that of dialectical laws defined by him only through examples. Thus the "dialectical law of transformation of quantity into quality" was valid for Engels as an empirical physical law (such as the assertion "water at the tempeature of 100° Celsius boils" - because the experiment repeated a large number of times brings the same result) as well as an economical law asserting that every sum of money cannot constitute a capital because a minimal quantity of funding is needed to exert a sufficient financial pressure. For Engels, the principles of materialism must be not the departure point of a scientific inquiry but its final result. The methodology was not to apply such principles to nature or human history but on the contrary to infer the principles from the study of nature or human history (Engels, 1959; p.54). The need for a dialectical scientific method comes from the fact that it was usual for classical science to observe

> *natural objects and processes in isolation, apart from their connection with the vast whole; of observing them in repose, not in motion; as constants, not as essentially variables; in their death, not in their life.* (Engels, 1959; p.34).

In the word "dialectics", Engels saw the principle for any motion of nature, history or thought, summarized in three basic laws : the "transition from quantity into quality", the "law of mutual interpenetration of opposites", the "law of the negation of the negation". Graham (Graham, 1987; p.29-30) insists on the fact that Soviet philosophers have often, voluntarily or not, neglected the continuous revision of the methods and results of a scientific approach implied by Marx's conception of science. This produced a dogmatism which was used to support political and bureaucratic interests and squashed the Soviet way of thinking for decades. It is remarkable in this respect that Lenin's *Philosophical Notebooks* (Lenin, 1961) contain many notes written by Lenin in order to improve and consolidate his own philosophical approach of these questions after the publication of his book *Materialism and Empiriocriticism* (which became the basis for the Soviet dogma in the 1920s) and remained rather underestimated for long in USSR and translated into English only in 1961. These notes show a Lenin much less dogmatic than the one drawn by *Materialism and Empiriocriticism*, admitting for instance a place for "fantasy even in the most exact sciences" (see (Selsam and Martel, 1963)).

2- The difficult relationship between the Intelligentsia and the new power

The launching of the project of the LSE in 1924 happened in a particular context whose instability simultaneously favoured the genesis of the enterprise and required a great flexibility from the editors and collaborators in order to adapt to the environment. We shall later show that this had consequences even for the entries about mathematics.

The relationship between the Bolshevik regime and the Intelligentsia had in fact been complicated since the very conquest of the power in 1917. To the eyes of the most sectarian Bolsheviks, it was a typical conflict of classes as the Intelligentsia, in its bourgeois way of living and thinking, proved to be a product of the old czarist society. One of the most crucial aspects of this opposition was the total defiance for letting the Intelligentsia dealing with the education of youth in an adequate spirit. The harsh period of war communism met with the peak of these tensions with the brutal decision for a proletarization of the whole scientific and technical personnel.

The Agitprop (*Агитационно-пропагандистский отдел* - Bureau for agitation and propaganda) was founded in 1920 under the supervision of the secretariate of the central committee in order to "organize, unite and direct all the oral or written work of propaganda and agitation" of the party[4], and this political propaganda was highly concerned with educational questions. It was especially efficient to denigrate the "old" and "bourgeois" specialists in the institutions of education, created in parallel to the old institutes and universities after the revolution, and destined to the education of "red" specialists and to the proletarization of universities : the Socialist Academy (*социалистическая академия*) created in 1919 which became, at the end of 1923, the Communist Academy (*коммунистическая академия*), or the workers' faculties (*Рабочие факультеты*) after 1920[5]. Moreover, this institutional politics was often accompanied by political violence. There were many press campains and show-trials with members of the Intelligentsia as targets. The GPU (*Государственное политическое управление* - State political direction, the State police) established a strict surveillance of scientific technicians who were easily accused of sabotage.

The period saw a drastic silence imposed to academic specialists considered as bourgeois representatives who damaged socialist edification. The targets were in the first place specialists in humanities : historians, economists or philosophers judged irretrievable by the Bolsheviks. A famous example is the "boat of philosophers" in 1922 on which many academics left Russia for a lifelong exile in Western Europe, such as Nikolaï Berdiaev or Lev Shestov. As Trotsky declared, "there was no sufficient pretext to shoot them, but it was no longer acceptable to bear them" (quote from (Ossorgin, 1955; p.183)). In his 2007 thesis, I.Kazanin emphasizes that one of the most efficient means of pressure used to transform intellectuals into pariahs who depended on the regime good-will, was to forbid their children to study in universities or institutes.[6]

Due to that ideological politics, the country was worn out in 1921 and Lenin decided the radical change of the New Economical Policy, to which we shall come back at length later. The incompetence of new economical decision-makers issued from the proletarization of economy was often such, that it produced a total destructuration of the means of production, already considerably weakened by the Great War and the Civil war. It was observed for instance that a lot of highly qualified engineers remained unemployed in Moscow or Petrograd, and nothing was organized to hire them for the huge construction sites in the inner country, such as in Siberia for instance. Harsh disputes broke out at the top of the Bolshevik party between the "inflexible communists" (*жёсткие коммунисты*), as Jozef Unszlicht who asked for an always stronger control of the Intelligentsia by the GPU, and the "liberal communists" who accused the first ones to lead the country to ruin and not to take into account Lenin's advice that " the best organizers and the top experts [could] be utilised by the state [also] in the old way, in the bourgeois way (i.e., for high salaries)" (Lenin, 1918).

---

4  Оргбюро ЦК 27 ноября 1921 г

5  See (David-Fox, 1997) and (David-Fox, 1998).

6  See in particular (Kazanin, 2007; p.325 *et seq.*).

In December 1921, Vladimir Vasilievich Oldenborger, a 58 years-old highly qualified hydraulic engineer, described by Lenin himself as the "people's commissar for water", was driven into suicide because a smear campaign organized by some local responsibles of the Bolshevik party who accused him of counter-revolutionnary sabotage. This event created a real shock in the party and was commented on by Lenin himself ("a shocking affair : we must sound the alarm" - see (Bailes, 1978; p.60-61)) and in many press articles, in particular in a momentous article of the Pravda on 3 January 1922 claiming for the immediate termination of such wastes[7].

As another example, let us mention that exactly at the same time, exasperated by numerous complains about how the party local organizations hindered the work of the statisticians in many provinces by depriving them of decent office buildings, the president of the Executive central commitee, Mikhail I. Kalinin, called to order the local responsibles on 21 January 1922

> *Due to the execution of important tasks by the statistical offices established by the decrees of the Central Executive Committee, of the Supreme Soviet of the National Economy and of the Council of labour and defense, and in connection with the conduct of the New Economic Policy, the Presidium of the central Executive Committee draws attention on the inadmissibility of these expulsions, transfer and crowding of the offices mentioned to other premises without their agreement.*[8]

3 - The NEP turning point

Lenin's decision to launch the NEP led to a partial and complex return to free-market economy from 1921 which attained its heights around 1925. This radical change allowed the liberal communists to decree a whole series of reforms including the relaxation of the politics of "class selection" and this led to a progressive normalization of the situation of the Intelligentsia, and especially of engineers and scientists. As Kazanin mentions

> *The agenda of many meetings of the Political Bureau and their decisions during the years 1924 - 1925 suggest that the leadership of the country was seriously concerned with the problem of reconciling the interests of the government and the technical intelligentsia, because its production efficiency, to a large extent ensured the economic and political stabilization in the country.* ((Kazanin, 2007), p.343)

It appeared of vital importance for the regime to bring the technical qualification of the specialists at the forefront and to let the question of their strictly political and social orthodoxy in the backgound, at least for a while.

> *During the meeting of the Political Bureau on 11 December 1924, was mentioned that state apparatus in its activity has to use experts not only from*

---

[7] See (Kazanin, 2007; p.316, especially note 722) for details and references about the reactions to Oldengorger's suicide.

[8] Quoted in (Mespoulet, 2001; p.240).

> *party workers -- the number of which in Soviet institutions is insignificant -- or from non-party workers, but, above all, from intellectuals and elements belonging to other classes, even if they are often alien to us, without which, the state apparatus can currently not do.[9]*

A calming down of the tensions between the authorities and the intellectuals was looked for so that the old specialists would be in a position to prepare the future executives coming from among workers and peasants (see for instance (Kazanin, 2007; p.348 *et seq.*).

The years of the NEP were an occasion for great political pragmatism. In August 1925, a report sent to the Central committee of the Bolshevik party proposed a series of measures to settle down the suitable conditions for an harmonious collaboration between the new executives and the former specialists so that a transfer of experience would be guaranteed which "can be realised only through common practical work during a significant time under the supervision of the old specialists".[10]

II - **The encyclopedic project of the Large Soviet encyclopedia (LSE)**

1- The stock company "Soviet encyclopedia"

A significant measure promoted by the NEP was the reopening of private publishing houses in order to improve the publishing activity in USSR and to give new platforms to the representatives of the "old" intellectual class for which access to publications controlled by the political sphere, such as the journal "Under the banner of Marxism" (*Под знаменем Марксизма*) published by the Communist Academy, was difficult. It is worth observing in passing that one must nevertheless not think that these publications presented a uniform opinion about science in the 1920s. On the contrary, the Communist Academy experienced harsh debates, and the journal "Under the banner of Marxism" exposed various aspects of controversies, sometimes with a rather abrupt tone. In the paper (Mazliak and Perfettini, 2016) one of them is studied about the so-called Marxist vision of randomness (to which we shall come back later in the third section) between E.Kol'man and V.I.Orlov. In fact, both Kol'man and Orlov claimed themselves to be loyal supporters of the regime, and the journal was banned for non politically comitted academics, and obviously to opponents.

In 1921, the celebrated publisher I.D.Sytin (И.Д.Сытин - see (Ruud, 1990)) observed that the volume of publications released by private publishers had become insignificant in comparison with that of the Gosizdat (Государственное издательство РСФСР - State publishing house of the Russian federation) which had been created in 1919. However, Sytin commented that

> *by ordering from private publishers the realization of books on a contract base, the Gosizdat would give them an opportunity to have a stronger situation. And at the same time this would enable to make use of the rich experience of the old publishers.[11]*

---

[9] Quoted by (Kazanin, 2007), p.345

[10] See (Kazanin, 2007; p.351).

[11] Quoted in (Koriakin, 2011; p.47).

The man who was to become the editor in chief of the LSE, Otto Yulievich Schmidt, to which we shall return in the next section, also observed the deficiencies of the Gosizdat whose head he was since 1921.

> *The Gosizdat must prove that it publishes books easily, well and cheaply and then it will, of course, be out of reach of competitors. (...) In Russia takes place a very curious experiment. We build the largest publisher in the world, but we do not give any commercial aim to it but only cultural and political ones.*[12]

Another important decision inspired by the NEP was the reintroduction of stock companies (*акционерные общества*) as a form of entrepreneurship in order to facilitate the elaboration of financial packages to support the activity of the company. The 1920s experienced the creation of a lot of private and cooperative publishing houses and a drastic increase of the journals and books production (see for instance (Kuznetsov, 2006; especially 111-112)).

A stock company was founded under the name "stock company Soviet encyclopedia" (*акционерное общество Советская Энциклопедия*), under the direct responsibility of the Presidium of the Central committee of the party, with a series of stockholders including the Gosizdat, the cooperative publishing houses "Questions of labour" (*Вопросы труда*), "Worker of education" (*Работник просвещения*), the publishing house of the people's commissariate for the workers' and peasant's inspection, the publisher "News of the Central election commision of the Russian Federation" (*Известия ЦИК СССР*), the publishing house "Pravda" (*Правда*), the stock company "International book" (*Международная книга*), the national bank of USSR, the Commercial and Industrial Bank of USSR (*Торгово-промышленный банк СССР*).

The history of the stock company "International Books" (*Международная книга*) for instance is representative of the urgence felt by Soviet leaders at the time of the NEP for improving the access to foreign litterature after years of complete isolation during the war communism. On 14 June 1921 V.I. Lenin signed the decree "Law on Acquiring and Distributing Foreign Literature", on the basis of which a Russian-German joint enterprise called "Book" (*Книга*) was established in Berlin. Its goals were import and export of books and other printed materials. Practically, the capital belonged entirely to the Soviet Union, and its head was USSR's trade representative to Germany, B.S. Stomonyakov. In 1922 was opened a branch in Moscow under the name "International Book", which on 11 April 1923 by special decree of the USSR Council of Labor and Defence was reorganised into a joint-stock corporation with the same name. All shares of the newly established corporation were again owned exclusively by the government of the Soviet Union, thus making it an entity with independent commercial agenda under the control of the state, a typical ambiguity of the NEP period. A list of merchandise started to include books, records, audio and video technology, machine parts, antiques, precious stones, philately, collectable coins and banknotes.

As mentioned above, the access to foreign technical litterature was a constant worry for Soviet leaders during the NEP. One reads in the minutes of the meeting of the Politburo from 13 August 1925

---

[12] Quoted in (Koriakin, 2011; p.47).

> *For this aim, an assistance was given by some trusts, enterprises or institutions, for publishing specialized journals and technical litterature, for giving access to foreign editions to specialists and organizations by buying subscriptions, so that they can be used by the technical personnel in order to facilitate the creation of necessary bonds with foreign scientists.*[13]

2- Otto Yulievitch Schmidt

As mentioned before, it was the former director of the Gosizdat, Otto Yulievitch Schmidt, who was put at the head of the encyclopaedic project. The aim of this section is to give some information about this amazing personality. Further details can be found in the biography (Koriakin, 2011) and in (Gliko, 2011).

Born in Belorussia in 1891 in a family of German descent, Schmidt studied mathematics in Kiev's university and began there his brilliant scientific career with some profound research in group theory, following his study of Jordan's treaty on substitutions. After discovering Remak's theorems on the decomposition of finite groups, Schmidt proposed several extensions of these results, and proved in particular a remarkable theorem[14] found at the same time by Krull and named afterwards after them both (*Schmidt-Krull theorem*) - see (Hungerford, 2008). Schmidt published his achievements in the book *Abstract group theory* published in Kiev in 1916 (Schmidt, 2015).

At the same time, he became interested in political action for the improvement of university conditions in the first place, but after the February revolution, with a more general concern about the future of Russia. In Summer 1917, he decided to come to Petrograd and to work for the Provisory Government on the question of food supplies. After the Bolshevik coup, he succeeded in proposing his services to the new people's comissariat for food supplies, probably a sign of his skillfulness at this position, sufficiently remarkable to have made him unavoidable. The opinion held about Schmidt at the top of the regime seems to have remained excellent in the subsequent years as it was at the request of Lenin himself that Schmidt was put at the head of the Gosizdat. It is not obvious that Schmidt scientific education played a role in this appointment, though Schmidt himself seemed convinced of the fact. In 1922, he wrote in a booklet published for the fifth anniversary of the Gosizdat (Schmidt, 1922)

> The upheaval and character of the extremely wide turn we experienced are reflected particularly in the scientific literature. We shall build the socialism on a scientific basis, on the basis of the Marxist theory and on the Marxist transformation of all the great discoveries of science.[15]

---

[13] Quoted in (Kazanin, 2007; p.350).

[14] The present paper is not the right place to comment on this result. Let me only mention briefly that it says that If G is a group that satisfies ascending and descending chain conditions on normal subgroups, there is a unique way of writing G as a direct product $G_1 \times G_2 \times \ldots \times G_k$ of finitely many indecomposable subgroups of G.

[15] Quoted in (Koriakin, 2011; p.48).

On 17 April 1924, the Central committee of the party approved the project of publication of the Soviet encyclopedia but the final choice of Schmidt as editor in chief would take place only on 15 January 1925. Maybe a cause of this long time for decision should be looked for in Schmidt's troubles at the head of the Gosizdat where he was opposed to supporters of a strict propagandistic aim for the publications of the house. In November 1924 for instance[16] Schmidt openly complained that some people seemed to wish the Gosizdat continue working as during the period of war communism, to the detriment of its financial situation. Ten days later, he was dismissed from his post by the *Narkompros* (People's commissar for education) and this may have facilitated his appointment at the head of the private company "Soviet encyclopedia".

We shall briefly describe the first editorial board of the encyclopedia in the next subsection. We shall in particular see that a large majority of its members were victims of the political storms experienced by Soviet Union in the 1930s. It is therefore slightly surprising that Otto Schmidt could remain at the head of the enterprise almost until the end (he resigned in fact in 1941), despite his proximity with Bukharin and even, to a certain extent, with Trotsky. Maybe Stalin thought it was useless for the regime to touch an internationally too well-known scientist. But above all, Schmidt himself had the wisdom, as soon as the end of the 1920s, not only to make a brilliant come back to mathematics (he was appointed to the newly created Chair of higher algebra at Moscow university in 1929 and remained there until 1949), but also to participate to long-distance scientific exploratory expeditions such as the German-Soviet expedition to the Pamir (1928) and afterwards the long expedition in the Arctic (1930-1934), which maintained him far from the internal struggles tearing the party apart at the turn of the 1930s.

Schmidt died in 1956 in Moscow, a rare example in Soviet Union of a person who was close to the power circles between the 1920s and the 1940s, without having been repressed.

3- The first editorial board

When the Soviet encyclopedia was launched, the editorial board headed by Otto Schmidt comprised thirteen members including :
N.I.Bukharin(*Н.И.Бухарин*),V.V.Kuibyshev(*В.В.Куйбышев*), M.N.Pokrovskij(*М. Н. Покровский*), G.I.Brojdo(*Г. И. Бройдо*), N.L.Mechtcheriakov(*Н. Л. Мещеряков*), L.N.Kritzman(*Л.Н.Крицман*),Yu.Larin(*Ю.Ларин*),G.M.Krzhizhanovskij(*Г.М.Кржижановский*),V.N.Miliutin(*В.П.Милютин*), N.Osinskij(*Н.Осинский*), E.A.Preobrazhenskij(*Е. А. Преображенский*), K.Radek(*К. Радек*), I.Stepanov-Skvortsov(*И. Степанов-Скворцов*).

Though it is a little digression from the central topic of the present paper, I think relevant to give some information on each of these members so that the reader can have a better insight on the spirit in which the first edition was launched. The following list intends to provide only some very limited facts about the lives of the members of the board, mostly in the years considered in this paper - namely the 1920s. The interested reader will easily find complements elsewhere if desired.

One observes a relative ideologic homogeneity in the board, as all the members in it claimed to be orthodox Marxists, and in general were persons having a political activity at the top of the State. The LSE was elaborated under the supervision of members of the "old guard" of the Bolshevik party, often even close acquaintances of Lenin.

---
[16] See (Koriakin, 2011; p.49).

It is remarkable that few members of the editorial board were academics, though they were all intellectuals (sometimes self-made intellectuals). The interest for economy and technique, emphasized in the preface as mentioned earlier, is visible in the choice of the participants.

A significant sign of the violence of the Soviet politics in the 1930s, it is seen that out of the 14 members of this committee, seven were eliminated in the 1930s (in general shot during the repression years 1937-1938), a proportion made all the more dramatic by the decease of three members before 1932, that is before the worst Stalinist repressions. Moreover, one member of the board, G.I.Brojdo, was condemned and sent to a camp, but had the exceptional fortune of coming back after Stalin's death. The repressed members of the board are highlighted in the following list by a "[r]" following their name.

- Nikolai Ivanovitch Bukharin (1888 - 1938) [r]

*One of the most important historical leaders of the Bolshevik party in the 1910s and the 1920s, Bukharin had met Lenin in 1912 and they became close friends since then. He was deeply interested in economy, had an excellent knowledge of Marxist writings and was considered in the 1920s as the intellectual endorsement in the party. Bukharin was a prolific author on revolutionary theory. He made several notable contributions to Marxist–Leninist thought, most notably The Economics of the Transition Period (1920). He was a founding member of the Soviet Academy of Arts and Sciences, and a keen botanist. His primary contributions to economics were his critique of marginal utility theory, his analysis of imperialism, and his writings on the transition to communism in the Soviet Union.*
*At the beginning of the 1920s, he was a go-between between Lenin and Trotsky. He became the foremost supporter of the New Economic Policy, to which he was to tie his political fortunes.*
*After Lenin's death, he was seen as one of his successors and became close to Stalin whom he supported in his fights against Trotsky, Kamenev and Zinoviev. It was Bukharin who formulated the thesis of "Socialism in One Country" put forth by Stalin in 1924, which argued that socialism (in Marxist theory, the transitional stage from capitalism to communism) could be developed in a single country, even one as underdeveloped as Russia.*

- Valerian Vladimirovitch Kuibyshev (1888-1935) [r]

*A specialist in electrotechnique, Kuibyshev was very early engaged in political agitation with the social-democracy and then with the Bolsheviks (he had in particular organized an illegal printing house in his town Kuznesk in Siberia). After the revolution, he was from 1921 in charge with the practical aspects of the electrification of the country as the head of the Glavelektro (the direction for electrification. Between 1923 and 1926 he was the People's commissar for the workers' and peasants' inspectorate.*

- *Mihhail Nikolaevitch Pokrovskij (1868-1932)*

*An historian from Moscow, Pokrovskii had not been allowed by the tasrist regime to defend his thesis considered as subversive. He became close to Lenin quite early in 1906 and worked for the Bolshevik paper "The worker". He left Russia in 1908 and came back only in 1917 after the February revolution. He was since 1918 one of the main organizers of the Socialist (Communist after 1923) Academy, and was the editor in chief of several newly founded journals as "Red archive", "Marxist historian" and "the class war". He was often an official representative of Soviet historical school in international conferences as in Oslo in 1928, the first time Soviet Union was officially invited to a congress.*

- Grigori Isaakovitch Brojdo (1883-1956) [r]

*Of Jewish ascendence, Brojdo participated to the 1905 Revolution, and to the February revolution in Tashkent where he had been sent. He was a Menshevik until 1917 and became then a Bolshevik. In 1921 he became the first rector of the communist university for the eastern workers, was the adjoint to the People's commissariat for the nationalities of the Russian federation headed by Stalin and became one of the responsibles of the Gosizdat.*

- Nikolai Leonidovitch Mechtcheriakov (1865-1942)

*Very early engaged with the social-democracy and then with the Bolsheviks, he was the editor in chief of the journal "The Krasnoyarskij worker" until 1917 and then came to Moscow. After the October revolution, he was a member of the editorial board of the journal of the Moscow soviet and then of the "Pravda". He became since 1920 one of the responsibles for the People's commissariat for Education and then the editor-in-chief of the Gosizdat headed by Schmidt. Mechtcheriakov was at the origin of the organization of ideological censorship for publications. He was nominated in 1921 professor at the faculty of social sciences. He participated also to the Komintern meetings.*

- Lev Natanovitch Kritzman (1890-1938) [r]

*Kritzman was one of the main economists of the Bolshevik party but he was in the first place a chemist graduated from Zurich university before defending a thesis in economy. He had fled to Zurich after his participation to the 1905 Revolution. After the Bolshevik coup, he remained close to the top of the regime. He became a member of the presidium of the Gosplan in 1921, a member of the editorial board of "Pravda", an organizer of the Communist Academy. In 1928, he became the adjoint of the head of the Central direction of Statistics and the same year a member of the Commission for the study of the agrarian revolution in Russia.*

- Yuri Larin (1882-1932) (pseudonym of Mikhail Zalmanovich Lurie)

*Larin was a Menshevik until the end of 1917 and then became a Bolshevik. He was one of the organizers of the Gosplan in 1921. He was a fervent supporter of the suppression of money in USSR, with money replaced by a direct distribution by the State of goods and services. In 1923, he was the initiator of the project of creation of Jewish autonomous regions in Crimea, Ukraine and Bielorussia.*

- Gleb Maximilianovitch Krzhizhanovskij (1872-1959)

*Krzhijanovskij was a specialist in electrotechnique and during the 1910s directed the construction of electricity plants around Moscow and in Saratov. He had also been engaged in the revolutionnary movement since his youth and was in close relationship with Lenin since 1893. After the October Revolution, he was in charge with the development of energy distribution in Russia and had on 26 December 1919 a famous exchange with Lenin about the electrification of the country from where Lenin drew his famous slogan "Communism is the soviets and the electrification of the whole country". In 1921, he became the first president of the Gosplan and as such supported the idea of concentrating the production and the distribution of energy in large production plants and regional energetic networks.*
*Between 1923 and 1926 he was the rector of the Moscow institute of mechanics. When Stalin decided to stop the NEP for a centralized planification, Krzhizhanovskij opposed this line and suggested a decentralized planification. It is worth noting that at the end of the 1930s he had the luck not to be arrested but only removed from his activities and kept under strict supervision.*

- Vladimir Pavlovitch Miliutin (1884-1937) [r]

*Miliutin was also a revolutionnary of the first hour who became one of the top leaders of the Bolshevik party after the Revolution. He was member of the economical commission for the North-West zone in 1921, the representative of the Komintern in Austria and in the Balkans in 1922, a member of the people's commissariat for Workers' and Peasants' inspectorate between 1924 and 1928 and he worked for the Gosplan after 1928.*

- N.Osinskij (*pseudonym of* Valerian Valerianovich Obolenskij) (1887-1938) [r]

*In 1917, after the Bolshevik coup, the revolutionnary activist Osinskij was mandated to crush the "sabotage" organized by the employees of the Russian State Bank and was nominated in December as head of the new State Bank of Soviet Russia. In 1920, he became a member of the People's commissariat for food supply. In 1923 he was a plenipotential representative of Soviet Union in Sweden. Then he worked for the Gosplan after 1925 and from 1926 was the head of the Central direction of statistics.*

- Evgenii Alexeevitch Preobrazhenskij (1886-1937) [r]

*In 1918 Preobrazhensky joined the Left Communists faction opposing the draconian peace with Germany established by the Treaty of Brest-Litovsk. It was at this time that he became closely affiliated with Nikolai Bukharin. In 1920–1921 he was Secretary of the Central Committee; in 1921 President of the Financial Committee and a member of the Council of People's Commissars of the Russian federation; Chief of the People's Commissariat of Education. Through the 1920s he was a leading Soviet Economist, developing the plan for industrialisation of the country and an opponent of the NEP.*

*He co-wrote the book The ABC of Communism with Nikolai Bukharin, with whom he would strongly disagree on the industrialization issue. He also wrote The New Economics, a polemical essay on the dynamics of an economy in transition to socialism, Anarchism and Communism and The Decline of Capitalism.*

*In 1924 he became one of the editors of the newspaper Pravda in 1924, a supporter of Trotsky as member of the Left Opposition. In the years 1924–1927 he was a member of the Board of People's Commissariat of Finance. After 1927, expelled from the party "for the organization of illegal anti-party printing house" and from January 1928, was sent to the Ural Mountains and worked in the planning agencies. In summer 1929, together with Karl Radek and Ivar Smilga Preobrazhenskij he wrote a letter claiming an "ideological and organizational break with Trotskyism". In January 1930, he was restored to the party and appointed to the Nizhny Novgorod Planning Committee and in 1932 became a member of the Board of the People's Commissariat of the Light Industry, and acting head of the People's Commissariat of State Farms.*

- K.Radek (1885-1939) [r]

*Radek was born in Lemberg, then in Austro-Hungary to a Jewish family. He became engaged in politics very young, and came to Russia to participate in revolutionnary movements. In December 1918, he participated in the discussions and conferences leading to the foundation of the Communist Party of Germany (KPD). On his return to Russia, Radek became the Secretary of the Komintern, taking the main responsibility for German issues. Radek was part of the Left Opposition from 1923, writing his famed article 'Leon Trotsky: Organizer of Victory' shortly after Lenin's stroke in January of that year. Later in the year at the Thirteenth Party Congress he was removed from the Central Committee. In the summer of 1925, Radek was appointed Provost of the newly established Sun Yat-Sen University in Moscow, where he collected information for the opposition from students about the situation in China and cautiously began to challenge the official Komintern policy. Radek was sacked from his post at Sun Yat-Sen University in May 1927 and expelled from the Party in 1927 after helping to organise an independent demonstration on the 10th anniversary of the October Revolution with Grigory Zinoviev in Leningrad. In early 1928, Radek was deported to Tobolsk and then to Tomsk. On 10 July 1929, alongside other oppositionists Ivar Smilga and Yevgeni Preobrazhensky, he signed a document capitulating to Stalin, was re-admitted in the party in 1930 and was one of the few former oppositionists to retain a prominent place.*

- Ivan Ivanovitch Stepanov-Skvortsov (1870-1928)

*Skvortsov-Stepanov was one of the oldest participants in the Russian revolutionary movement and a Marxist writer who became a Bolshevik in the winter of 1904. When the journal "the struggle" (Борьба) was published in November 1905, Skvortsov-Stepanov was a member of the editorial board. In 1906 he was a delegate to the Fourth Congress of the Russian Social Democratic Labour Party, where he supported Lenin. During the period 1907–10, he favoured the Mezhraiontsy faction, but later fell again under the influence of Lenin. He was repeatedly arrested and exiled for his revolutionary activities. Following the Revolution of 1917 he became the People's Commissar for Finance of the RSFSR. He was also, during the 1920s, one of the most ferocious organizers of the struggle against religion. Upon his premature death in October 1928 after he contracted a severe typhoid fever, Stepanov was commemorated by Stalin as a "staunch and steadfast Leninist".*

4- First tensions

As explained before, the launching of the LSE benefited from the narrow window of the NEP which did not last long. At the end of 1928, tensions gradually increased, and a violent campaign orchestrated by the regime against the "bougeois specialists" began. For a general overview on this question, consult for instance (Werth, 2012) (chapter 6) and (Krementsov, 1997). Aspects concerning mathematics have already been studied several times : see in particular (Tagliagambe, 2003) (in particular p.88), (Vucinich, 2000), or (Seneta, 2004). An enterprise as the LSE was hit with full force by this campaign.

Schmidt's choices for the publication and the composition of the editorial board placed the LSE in the front row to be under attack when began the idological crispation in the beginning of the 1930s. In 1929, the journal Natural science and Marxism (*Естествознание и Марксизм*) devoted to the study of natural sciences from a Marxist point of view was founded under the direction of Schmidt as an extension of the Soviet encyclopedia with more room available for debates. In 1931, Schmidt was replaced by E.Kol'man and the journal was refounded under the title For a Marxist-Leninist natural science (*За марксистко ленинское естествознание*) which in its first issue declared war to Schmidt. His choices for the LSE were for instance violently opposed in an article by A.A.Maksimov ((Maksimov, 1931) - see especially p.73), asserting that the "science section of the LSE must be considered as anti-marxist. (…) To the ideological character of some observed errors must be related the errors made in choosing the authors." And the logician S.Yanovskaya, who, for some times had become E.Kolman's companion as a guard dog of the ideological purity of Soviet mathematics, composed a long diatribe against the way in which mathematics were presented in the encyclopedia. She wrote

> *If we examine the contents of the section on mathematics in the Encyclopedia, we find that not only there is in it no trace of critics towards idealism, but on the contrary, as rightly pointed Comrade Maximov, one discovers there its exaltation.*
>     *For the idealist, a precise system of axioms allows to define real numbers, which could be defined as well in another way, namely by so-called genetic method corresponding in all and for all to the historical development, linked as it is to the progressive movement, to the extension of the concept of number and of the operations that are made on it. From this point of view, however, it is not the genetic orientation which is fundamental and decisive: everything is reversed and the system of axioms comes to be considered as something which, like Minerva, comes fully dressed out of the head of Jupiter, and thus surges in an already completed form out of the head of modern mathematics. It is truly from this system, henceforth claimed to be the demiurge of mathematics, that is awaited the stabilization of the process and of the limits of the successive genetic expansion of concept of number. Discussions between genetic and axiomatic methods still continue to shake the philosophy of bourgeois mathematics, but of this we find no trace in the entry [axiom], despite it is of vital importance for us to explain the subordinate role of axiomatic in science and to emphasize that axioms are in fact not an initial point but rather the result of an analysis started at an already rather high level of scientific development. (...) For all this, I retain, in conclusion, that the entry "Axiom" should be seen as idealistic and not Marxist.* (Yanovskaya, 1931; p.79-81).

5- A periodization related to Soviet inner politics

The printing of the first edition of the LSE was set at 60,000 exemplaries. As already mentioned, the editorial board had chosen for the series the form of an encyclopedic dictionnary and a publication of the volumes in alphabetical order. This order was almost respected all along the 25 years period of publication. This organization produced unexpected difficulties to deal with in the meanders of the ideological control by the party during the 1930s. Some people for instance were deprived of a decent entry in the LSE such as D.Egorov, politically suspect and who would die in exile in 1931 (see (Seneta, 2004)). The very short article on Egorov (published in 1932, one year after the mathematician's death - he probably would not have had any entry at all if he were still alive) is deeply derogatory, mentioning that Egorov's mathematical achievements were not important, and that he was above all a representative of the reactionnary Moscow mathematical school (we shall come back to this school in the third part when commenting on Bugaev).

There are also narratives dating from these years explaining how libraries and subscriptors received isolated pages to be pasted on the original pages of a volume of the LSE containing a politically banned text. A famous example (concerning however the second edition, I was not able to discover a testimony of the kind for the first edition) is Lavrentii Beria, as reported in Sharansky's memoire (Sharansky, 1988).

> *A blatant example of presentism appears in Natan Sharansky's description of modifications to the Great Soviet Encyclopedia. The first target of revision was Beria, a chief of the secret police who was executed for being a British spy. Subscribers to the encyclopedia were instructed to destroy the article on Beria and were provided with additional information on the Bering Strait to fill the gap in the pages. According to Sharansky, subscribers frequently received such missives.* (Starzyk, Blatz and Ross; p.465)

Even after the Stalin era ended, facing such a massive information as an encyclopedia contains remained a tough difficulty for a regime wanting a strict ideological control on the publications. The historian of book edition F.N.Petrov (Petrov, 1960) for instance, even as late as in 1960, was careful when he examined who should be advised to read the first edition of the LSE. He wrote

> *The first edition of the LSE keeps its importance up to the present time for providing information. In its biographical and historical aspect it can serve as a source of information for historians and researchers. But one cannot recommend it for large circles of readers as much of the contents include ideological and political errors and do require changes.* (Petrov, 1960; p.136)

Indeed, for Petrov, ideological aspects were inherent to the nature of the encyclopedic project. He wrote that the LSE "contains assertions with a political and international character. Therefore each formulation must be closely verified and must be perfectly conform to the ideological and political problems of our country" (Petrov, 1960; p.134). And he openly claimed that "the reader wants to receive methodological or political instructions (*установки*) based on Marxist-Leninist teaching in order to clarify events and fact occuring in nature and society" (Petrov, 1960; p.135).

Hence, the alphabetical order chosen for the publication of the volumes makes the analysis of the volumes of the first years all the more significant to perceive some ideological and scientific debates of the period 1925-1930, when a (very) relative freedom of speech still left some room for them. It is clearly not in mathematics that the penetration of the ideological debates was the most obvious, but in a way this is also why it is so appealing to drive the more or less subtle political instillations out in the entries about mathematics. And it would certainly be a mistake to think that mathematics were spared altogether. As Hayek commented in 1944

> *Totalitarian control of opinion extends, however, also to subjects which at first seem to have no political significance. Sometimes it is difficult to explain why particular doctrines should be officially proscribed or why others should be encouraged, and it is curious that these likes and dislikes are apparently somewhat similar in the different totalitarian systems. In particular, they all seem to have in common an intense dislike of the more abstract forms of thought-a dislike characteristically also shown by many of the collectivists among our scientists. Whether the theory of relativity is represented as a "semitic attack on the foundation of Christian and Nordic physics" or opposed because it is "in conflict with dialectical materialism and Marxist dogma" comes very much to the same thing. Nor does it make much difference whether certain theorems of mathematical statistics are attacked because they "form part of the class struggle on the ideological frontier and are a product of the historical role of mathematics as the servant of the bourgeoisie", or whether the whole subject is condemned because "it provides no guarantee that it will serve the interest of the people". It seems that pure mathematics is no less a victim and that even the holding of particular views about the nature of continuity can be ascribed to "bourgeois prejudices". According to the Webbs the Journal for Marxist-Leninist Natural Sciences has the following slogans: "We stand for Party in Mathematics. We stand for the purity of Marxist-Leninist theory in surgery."* (Hayek, 1944; p.165-166)

## III - Some aspects of the first mathematical entries in the LSE

1- V.F.Kagan and the natural science section

The first editor in charge of the part "natural and exact sciences" (естествознание и точные науки) was Venyamin F. Kagan. Born in Lithuania in 1869, Kagan was a rather typical member of the Jewish intelligentsia. Involved in the democratic motion of students at the university of Saint-Petersburgh (he was expelled from university for that purpose in 1889), he obtained his graduation in mathematics under the direction of Markov and Posse. He was then appointed as professor of mathematics in Odessa in 1897 and remained there until 1923 when he was given the Chair of Differential geometry at the University of Moscow. He was at the head of the Journal of experimental physics and elementary mathematics from 1902 until 1917. He had also an important publishing activity at the head of the house Mathesis specialized in the printing and diffusion of mathematical texts, described in (Lopschitz and Rashevskii, 1969) as the most important in Russia (see also (Rikun, 2012)). For this reason, Otto Schmidt proposed to Kagan a collaboration with the Gosizdat. He sent the following rather pressing letter to his colleague from Odessa in 1922

> *Dear Veniamin Fedorovich,*
> *There is no doubt that the Mathesis edition has offered to the country the greatest achievements in the field of books of exact sciences. The extension of this work is an essential cultural requirement.*
> *On the other hand, the Gosizdat feels necessary to offer scientific and*

> *vulgarization litterature. Want it or not, we face the obligation of publishing similar books and even sometimes to republish the same as you. I fear it does no good. The Gosizdat, being an overwhelming economic strength, will annihilate the revival of the publishing house Mathesis without using its know-how and traditions.*
> *That is why I propose the following. Why not merging you and us ? The Gosizdat would provide the capital for the revival of Mathesis under your leadership. This would constitute an autonomous section of the Gosizdat of the Russian Socialist Soviet Federative Republic. A similar experiment was conducted with the house "World Literature" (Gorki, Tikhonov) and gave great results. Both parties were fully satisfied.*
> *Please think about this and send me your opinion.[17]*

After he was appointed as editor in chief, Schmidt proposed Kagan to head the section of sciences of the LSE.

2- The mathematics of randomness in USSR in the 1920s

The end of the paper is devoted to examine four entries published in the first volumes of the LSE. Each of the four is in some way connected to the important debates existing in Soviet science in the 1920s about the status of the calculus of probability and its use in the scientific approach of phenomena. The discussions about the right place to give to the right place to give to randomness, and to its scientific measure, in the communist society under construction were indeed a great theme of reflections and exchange during the 1920s. Is is worth noting that the debate took place even in the ranks of supporters of the regime. This is for instance illustrated in the paper (Mazliak and Perfettini, 2016) where the harsch exchanges between E.Kol'man and V.I.Orlov are examined, after the publication by the latter of his volume on logic of natural sciences (Orlov, 1925), and the contestation by Kol'man of the Marxist orthodoxy of Orlov's views. The economical primacy resulting from the Marxist social conception led indeed to ask what margin of randomness was left politically admissible when the means of production were supposed to be under the absolute control of the State. In the economical sciences, any excursion beyond strict deterministic models was considered with a priori suspicion and essentially related to the existence of a market where private actors could speculate. It was for sure an originality of the Soviet scientific scene of the 1920s and 1930s that the mathematics of randomness were at the same time a topic in which Soviet science obtained blatant successes and one regularly subject to hard critics. The well known and dramatic case of the debates about Darwinism and genetics in the 1930s was also related to these questions - see (Krementsov, 1997) about genetics and lissenkoism. In (Mespoulet, 2001), Martine Mespoulet gives a detailed picture of the defeat of statisticians at the end of the 1920s to impose their methodology in front of the party organizations in search of figures of production matching the political agenda of the planification.

The aforementioned discussions between Orlov et Kol'man in the journal *Under the banner of Marxism* studied in (Mazliak and Perfettini, 2016) are another illustration of how the question was considered crucial during the 1920s. In these years, Soviet mathematicians felt necessary to prove that their mathematics were not "empty". The debate was also alimented by reflections about axiomatization. A.Ya.Khinchin for instance wrote the article (Khinchin, 1926) in *Under the banner of Marxism* to emphasize the importance of this "battle for the object" in modern mathematics. Khinchin explains how Weyl and Brouwer, when they wanted to "pitilessly expell everything which hides its emptiness under the veil of a perfect logical outside from mathematics" (Khinchin, 1926; p.184) did not wish to prove how some contemporaneous mathematics were pointless but wanted to show "a deep inner illness" of contemporaneous mathematics. In the Soviet society under construction, formalism was beginning to be considered with high suspicion (for details on these questions, see (Verburgt, 2016)).

---
[17] Quoted in (Lopschitz and Rashevskii, 1969; p. 19).

This may partly explain why at the time of the Stalinist turn of the 1930s, even a star mathematician as Kolmogorov felt necessary to make rhetorical efforts to convince his readers that, though deeply involved in probability theory, he was acting as a mathematician and was concerned only with the mathematical aspects, leaving to others the question of interpretation, connection to the real world and practical application of his research. He wrote for instance in his fundamental paper on analytical approach for Markov processes

> *It should be noted that the possibility to apply the schemes of deterministic or stochastically defined processes for dealing with any real process, has no connection with the question of whether this actual process was itself deterministic or random.* (Kolmogorov, 1931; p.3)

But, in the first years of USSR, scientists did not hesitate to expose various opinions on these questions, even if they claimed themselves to be faithful to Marxism. They used a variety of arguments to prove that they were on line with the Marxist-Leninist dogma, in particular to show that they were hostile to anything which could, from close of far, be called idealist. This was often a self-protection against a possible accusation of being too close to the activity of unreliable scientists such for instance those in Moscow with a religious background (see, among others, (Ford, 1991) and (Seneta, 2004)).

3- Comments on four entries

As said before, there were many debates on probability in young Soviet Union, and therefore entries in the LSE dealing with questions concerning randomness and its scientific estimation give a good insight on these matters. There is however another reason to examine the four articles we shall comment on. By the chance of the alphabetical order, they all four belong to the beginning of the alphabet and therefore were published around 1926 in the very first volumes of the encyclopedia. Two entries concern fundamental theoretical aspects : probability and large numbers (law of) and the two others deal with two mathematicians, one living and a foreigner, Emile Borel, the other dead and Russian, Nikolaï Vassilievich Bugaïev. The aim of the following comments is above all to illustrate the kind of balance the authors had kept between their scientific freedom and their necessary adaptation to the circumstances, before the harsh taking in hand by the Stalinist one-track thinking of the 1930s.

a- Probability (вероятность) *by A.Bowley and A.A.Khinchin*

The entry "probability" is an article in three parts : "mathematical foundations of the theory of probability", "calculus of probability", "application of the theory of probability". The second part, which contains more technicalities, was written by the young Aleksandr A.Khinchin who had recently begun his work on probability. He wrote a text without much originality. It contains four sections : Origins and development of the calculus of probability, Probabilities of compound and independent events, Probability of hypotheses and Bayes rules, Continuous probabilities (in which Khinchin presents a short exposition of Bertrand paradox to emphasize the necessity of a clear setting of the random experiment before any calculation).

More significant are the two other parts written by the British statistician Arthur Lyon Bowley. This is not a minor point as few foreigners, on the whole, were called to collaborate to the LSE. Bowley seems besides to have met a particular favour in these years in USSR. Though he was not a communist, and, up to my best knowledge, did not publicly express a whatsoever positive opinion on Bolshevism, a reason may have been that Bowley had been one of the initiators of statistics about the working class and its condition for living (see for instance (Bennett-Hurst and Bowley, 1915)). Also Bowley was a great supporter of using mathematics in economy (he published in particular a remarkable textbook of mathematical statistics for students in economy (Bowley, 1907)) : in a famous letter, his colleague and friend the economist Alfred Marshall teased him for this. He wrote : "I had a growing feeling in the later years of my work at the subject that a good mathematical theorem dealing with economic hypotheses was very unlikely to be good economics" and even advised to "burn the mathematics" as a rule for economists (Ekelund and Hebert, 1999; p.362). Bowley's opinion was exactly the contrary and the mathematical treatment for studying social problems was seen as a warrant against the subjectivity of a more discursive form. The interest for social problems, and the mathematical orientation of Bowley's works made probably him particularly acceptable in Soviet Union during the liberal times of the NEP. In absence of archivial material, it is difficult to know how Bowley was asked to participate to the LSE, who translated his contributions and if those were specially ordered for the LSE. His treatment of the two parts was besides not particularly engaged but one observes the accent he almost exclusively put on the frequentist approach. The application of the calculus of probability is legitimate only in the context of the law of large numbers. This approach, formalized in particular by von Mises at the beginning of the 1920s, was considered as the only one not contaminated by idealism. In his (Grundbegriffe, 1933), Kolmogorov would again carefully insist that he shared von Mises' frequentist point of view as the empirical justification of axioms which are only mathematical abstractions more easily handled with. He wrote

> *The reader interested only the purely mathematical development of the theory may not read this section (...). Here we limit ourselves to a simple exposition of the empirical origin of the axioms of the theory of probability and voluntarily leave the deep philosophical questions about the understanding of probability in the experimental world aside. For the exposition of the necessary hypotheses for applying the theory of probability to the world of real events, the author mostly follows the reflections of von Mises (...)* (Kolmogorov, 1933; p.3)

Let us finish this section by mentionning that the references given by Bowley include three modern Russian textbooks: (Markov, 1924), (Bernstein, 1927) but also (Lakhtin, 1924) by Leonid Kuz'mich Lakhtin who was a representative of the Moscow mathematical school who had just published his lectures. A natural question to ask is thus to understand why the board chose the foreigner Bowley for the entry instead of Lakhtin, for instance. The opposition against the approach of the Moscow school (see below the section about the entry "Bugaiev") may be a reasonable hypothesis. Apart from these Russian sources, the litterature includes (Poincaré, 1912), (Czuber, 1914), (Castelnuovo, 1919) and (Lévy, 1925) and Bowley's own book (Bowley, 1907).

b- Law of large numbers (больших чисел (закон)) *by M.N.Smit-Falker*

The second entry I would like to consider was written by the economist Maria N. Smit-Falker. She had been a convinced Bolshevik for a long time since 1907, and studied some years at London School of Economy before the Revolution. She met there Arthur Boley and remained afterwards in contact with him. This may be a hint to explain how Bowley happened to have written an entry for the LSE. Smit-Falker's tense relations with the world of statisticians in USSR had been studied in (Mespoulet, 2001; p. 292-293) : she was at the head of a department of the Supreme council for national economy between 1918 and 1920 and then became a professor at the Institute Plekhanov for national economy. Martine Mespoulet exposes how Smit-Flakner was convinced that the principles of rationalisation of the industrial production were transferable to the statistical activity. She wrote for instance :

> *For the numerous processes of recollection and treatment of data, rationalization of work and introduction of assembly-line work (коонверизация) must play a great role in order to reduce the waste in work. Up to now, we almost do not have any norm for productivity and no system of decomposition of statistical operations along a system of assembly-line.* (Smit-Falker, 1927; p.15-30)

The long article that Smit-Falker wrote for the LSE about the law of large numbers is oriented towards political economy both in the examples she presents and in the interpretation of the result. The following quotation shows how she wanted to emphasize that economy was the most essential application of the statistical method.

> *A collective is said to be statistical, if in it any character inherent to its members is unevenly distributed. For instance the value of the cultivated area of a peasant household enormously varies and for some properties is equal to zero (households without lands). The number of workers also varies from plant to plant.*

Moreover, Smit-Falker's text offers a striking example of how the rhetoric of excommunication of undesired people (here, members of the party), called by the sinister term of "purge" (чистка) which was so frequently used during the next decade and became a synonym of capital punishment, was already found in any kind of text, even an entry about a mathematical theorem in an encyclopedia.

> *In order to know the composition of a whole mass, it is necessaty to measure its totality or a sufficiently large part of it so that within this part will appear the connections which are present in the whole mass. For instance when a purge (чистка) of the party is decided the ratio between the number of members subject to exclusion and those not subject to exclusion in the individual cells can be very different. In some cells, one kind of party members prevails, in some other ones, another kind, and only by increasing the number of tested cells it is possible to refine the picture of the composition of the party as a whole.*

The statute of auxilliary of economy attributed to probability by the author is observed above all in her recurrent use of illustrations drawn from economic life. M.Smit-Falker wrote for

instance

> *The logical foundation of the law of large numbers is quite clear. On all members of a collective causes of a general character weight, while individual members can be affected by moments deflecting the impact of the general factors in either direction. For example the level of labour productivity of an entire set of factories and plants, generally reflects the overall level of development of productive forces and cultural skills of workers. But some factories can either keep up with the general level or, conversely, go ahead of it.*

Naturally, the Marxist dogma was conveyed whenever possible. Smit-Falker quotes for instance the following passage from Marx's capital[18]

> *Of the six small masters, one would therefore squeeze out more than the average rate of surplus-value, another less. The inequalities would be compensated for the society at large, but not for the individual masters. Thus the laws of the production of value are only fully realised for the individual producer, when he produces as a capitalist, and employs a number of workmen together, whose labour, by its collective nature, is at once stamped as average social labour*

This claim to Marx was besides not only seen as an ideological basis. It was also (and maybe above all) a weapon to disqualify classical - that is bougeois - statistics supposed to be at the service of capitalist oppression. The following quotation shows how Smit-Falker tried to oppose a kind of idealistic statistics supposed to be relevant only for hazard games in which conditions remain constant to the statistics needed in economy where such conditions do not happen.

> *In formulating the law of large numbers, the classical theory of statistics did not come from the observation of some social mass subject to change, but from the observations of cards, dice (gambling) or urns with black and white balls. It is possible to prove that the number of black and white balls in the urn is identical through a large number of repeated draws of balls from the urn (if the drawn ball is each time put back). (...) This formulation is quite correct at any time or place as the black and white balls, returned in the urn after the realization of the experiment, are not submitted to any influence from the outside; hence the resulting ratio is persistent in time. But in real, and especially in social collectives this can not take place. The composition is subject to continuous change in time. The party members do not come back to the party, as do black balls in the urn. (...)*

> *In biology until Darwin plant and animal species were considered fixed and established once and for all, political economy was dominated by eternal, immutable economic categories, and in statistics one was taught with the doctrine of so called "permanences". Economy proclaimed the imutability of such economic categories as capital and wage labor, statistics*

---

[18] Karl Marx, Capital, t.1, Chap.13.

> *proclaimed the stability of most statistical coefficients. The composition of social collectives seemed once and for all established and even considered as "the divine order of things" (Süssmilch). And in the work of some modern statisticians we still meet the doctrine of permanences, and for other statisticians the obtained relations are like "logical constants" or, in other words, relations, not resulting from a strictly limited experiment but applicable always and everywhere. Modern Marxist statistics, which always deal with the study of phenomena in the process of their formation, use the law of large numbers for the study of certain collectives in each given period of time, and there is no such a thing as a timeless effect of the law. Therefore, the relations obtained in collectives are not treated by modern statistical theory as a kind of "natural law" or as some "logical constant." Even in the case when the composition of a collective, observed at different time periods, is relatively stable, we are dealing only with slow change rather than with stability, and with an empirical constant rather than with a logical one.*

I shall now examine two biographical entries about two mathematicians who reflected about the evolution of probability theory at the turn 19th/20th century.

c- Emile Borel (Борель) *by N.Luzin*

We have seen that the works of the British statistician Arthur Bowley were quite positively welcomed in the young Soviet Union. It seems to have been also the case for the mathematician Emile Borel who benefited from a rather long entry in this first edition of the LSE, a noticeable fact for a living and non Soviet personality, as a strong pan-Russian (or at least pro-Soviet) tropism can be observed in the general economy of the publication. One may think that Borel's conceptions about the role of a mathematician in the city and his opposition to the most formalist aspects of the disciplin (Mazliak and Sage, 2014) made him a frequentable person. A sign of this favour is the series of translations of books by Borel published in the 1920s in Soviet Union, including for instance his well-known book *Le Hasard* (Randomness) (Borel, 1923) which was considered by Borel as a survey of his conception of the mathematical approach of randomness - see (Bustamante, Cléry and Mazliak, 2015). The editorial board asked Nikolaï Luzin for writing the entry about Borel. Luzin had followed Borel's lectures in Paris before the war and with his master Egorov and other students, was a passionate follower of the French mathematical works in theory of functions (Borel, Lebesgue…) on which they produced a series of important results during the 1910s. Luzin had founded in 1917 a seminar on these questions at the university of Moscow, the famous group Luzitania, where in the 1920s many future stars of the Soviet mathematics as Khinchin or Kolmogorov made their first steps.

In his article, Luzin underlines the role of Borel as being one of the first to understand the importance of Cantor's works but also to warn his fellow mathematicians against the risk of drift resulting from a purely logical approach. Luzin writes that Borel was positive about Cantor's theories

> *when these ideas were met with total disbelief. He first applied them for research on functions (Heine-Borel's theorem). However, with his inherent tendency to classical simplicity and concreteness, Borel warned scientists*

> *against their attraction for purely logical construction of infinite sets without an analysis of their relationship to reality. Borel's considerations («Illusion du transfini») were at the beginning not well understood, but further development of the theory of functions attracted the general attention on them.*

Luzin moreover mentions the wide selection of topics in which Borel was involved, with a special accent put on probability

> *Borel is keenly interested in many problems of mathematical physics, and in particular in the theory of probability, a field in which he is taking nowdays the edition of a series of monographs.*

Clearly, a realization as the *Treaty of probability and its application*, launched by Borel in 1922 (see (Bustamante, Cléry and Mazliak, 2015)) could not be suspected of any idealist tendency and Borel was thus quite acceptable within the Soviet scenery. Moreover, his political engagements, in particular at the Society of Nations with the Institute of Intellectual Cooperation (see for instance (Guieu, 1998)), would always make him very careful to maintain the contact with Soviet scientists, in particular when he was the head of the Institut Poincaré in Paris since 1928. He did his best during the 1930s to be able to invite Soviet mathematicians to Paris. He had a relative success at the beginning of the 1930s but remained helpless when Stalin decided to close the borders of the country.

d- Nikolaï Vassilievitch Bugaiev (Бугаев) *by V.F.Kagan*

With Borel, we had an entry about a living mathematician, foreign but considered as acceptable for the Soviet standards. With Bugaiev, we have an article about a dead mathematician, Russian but mostly unacceptable. And, in fact, it is probably this combination that he was Russian and dead which enabled to incude Bugaiev with a reasonable honesty in the LSE when his still living colleague Egorov, as mentioned above, began to be treated as a plague-stricken. Besides, it was V.Kagan himself who took charge of the text about Bugaiev and it may be why the necessary expression of hostility against him remains quite moderate while it may have been more violent if the article were written by a second fiddle wanting to give hints of submission to the regime.

Bugaiev was one of the founders of what has been called the philosophical-mathematical school of Moscow. Several of its members had deep connections with religious circles ; it is for instance in this school that a personality as Pavel Florensky would study at the beginning of the 20th century[19]. Bugaiev created a new discipline he called arithmology, a science of discontinuous functions aimed at a representation of the world richer than that of the Newtonian cosmology. Probabilities entered Bugaiev's system as an essential tool in order to go beyond arithmology. Here is for instance what Bugaiev declared in his conference of the first International Congress of Mathematicians in Zurich in 1897

> *Only a continuous education can restrict and weaken the limits of uncertainty in our judgments and actions. The very necessity of education that one can give to oneself already gets rid of fatalism of our theoretical*

---
[19] Numerous studies exist on P.A. Florensky, mostly in Russian and in Italian; see in particular (Betti, 2009).

> *views on man and nature. A certain element of randomness, appearing in our actions, introduces a contingency element in nature itself. Contingency thus comes on the stage, as an essential property of certain phenomena in the world. In the world there is not only the reign of certainty. Probability also has a place there. The doctrine of incidental phenomena or probability theory appears as an essential mathematical science in the general system of our knowledge. The philosopher must reckon with probability as much as with certainty. Probability theory must give answers when one cannot use analysis or arithmology, when we ignore the law of phenomena.* (Bugaiev, 1898; p.219-220)

Let us add that since the end of 19th century two approaches of probability coexisted in Russia. The one I mentioned in Moscow, and another one developed in Saint Petersburgh after Chebyshev by his disciples Markov and Lyapunov, more connected to application and distant from any metaphysical interpretation. Markov's violent hostility towards Bugaiev's favorite disciple Pavel Nekrasov and in particular towads Nekrasov's conceptions of probability, was part of the pre-revolutionnary Russian mathematical scene and this played some role in the modelling of an acceptable probabilistic theory in Soviet Union (see (Seneta, 2003)). One reads the following in Kagan's article

> *[Bugaev] believed that the doctrine of non continuous functions should constitute a great discipline, called by him "arithmology" which, he thought, was to cover the entire mathematical analysis, by taking over infinitesimal calculus. However, the studies performed by Bugaiev provided no reason for such broad generalizations. Meanwhile, Bugaiev put these views as foundation of his philosophical outlook, leading to the following. Determinism has its source in infinitesimal calculus. Laplace saw the justification of determinism in the existence of integrals of the differential equations of motion. But to Laplace and his followers was not known the arithmology which made explicit that in nature there are jumps, and contradicted the doctrine of determinists. On this basis Bugaiev exposed deeply metaphysical beliefs, and together with his students, of whom the most active was professor P.A.Nekrasov, created in Moscow a whole school of philosophy, with a clearly metaphysical direction, which had a great influence not only in mathematics, but also in the wider circles of Moscow scientists. Several representatives of the "school" conceived these philosophical deductions in relation with political views of a clearly reactionary nature. First rate Russian mathematicians as P.L.Chebyshev, N.A. Korkin and A.A.Markov were not inclined to this metaphysical constructions, and proved their inconsistency and did not recognize even any "arithmology".*

**Conclusion**

During the Central Committee meeting in July 1928, Stalin explained that the NEP was in a dead-end and that he was considering to require peasants to provide the efforts necessary to support a quick industrialization of the country. These words were sealing the rapid ending of the NEP but also the headlong rush towards the terrible years of the collectivization at the beginning of the 1930s. Bukharin, horrified by the perspective of terror and violence implied by Stalin's words attempted to resist by publishing on 30 September 1928 an article in the Pravda entitled "Notes of an economist" ; he tried to prove through a scientific analysis that the projects of creation of kolkhozes and of general planification of economy were extremely risky. Stalin obviously would not listen and decided to go forward. Bukharin had in fact signed his own death sentence. As soon as 1929, he began to lose gradually all his official positions until his complete isolation and elimination in the 1930s. He was one among hundreds of thousands academics and specialists who began to be repressed in the so-called "great turn" of the years 1928-1931. The nightmare had begun. In 1934, Stalin had eliminated all opposition, at least virtually before it became physically (see (Malia, 1999) or (Werth, 2012) among numerous books devoted to the history of the Stalinist period).

We have seen in our article how much the LSE was directly hit by these events. There were several complete changes of editorial board during the 1930s, A noticeable exception was Otto Schmidt who, apart from being safe during his long-distance exploratory expeditions, probably tried to be the more transparent he could as editor in chief of the LSE.

The years of the NEP appear thus to have been a parenthesis of relative quietness between the period of the war communism and the beginning of the Stalinist dictatorship. During these five or six years, there was some room left for academic debate, at least in the domains which were not the more exposed to a political interpretation. Even during the NEP, there was obviously no real freedom of speech on strictly political matters and certainly no hint at all for an opening of the political scene to any other party than the communist party. But on scientific questions, for instance, the new economical conditions, in particular for publishing books or for getting an access to foreign litterature, saw the emergence of vivid debates, sometimes even inside the structures created by the party.

This makes the publication of the volumes of the LSE during these years all the more interesting to understand the kind of thinking and, often, of utopia which reigned in Soviet Union during its first years about science. Moreover, the questions around randomness, as we have seen, acquired a special significance in the society without classes under construction because, through the prism of primacy of economy postulated by Marxism-Leninism, chance seemed to be a secondary concept if the state possessed all the control sticks of economical life. As all this happened at a time when the mathematical theory of probability was living a profound evolution, both because of problems on theoretical foundations and because of the creation of new concepts (processes, limit theorems), some original developments were proposed on the topic in USSR. This, certainly, is not completely unconnected with the prodigious rise of probabilistic studies during the 1930s, though, as we have already observed about Kolmogorov, with a special care for remaining far from any question of "concrete" interpretation.

In the present paper, we obviously have only studied a very specific aspect of the LSE. Such a gigantic source deserves much more work and it is certainly to hope that historians of science, and in particular of mathematics, will conduct other inquiries on many other aspects of the enterprise. The access to primary sources, such as the documentation about the launching of the enterprise (maybe the private correspondences of some members of the board, or texts about the project for instance) would be of capital interest. Finally, it seems that it was above all the Stalinist period which attracted more work about science in USSR. On one hand this is logical : the Stalinist dictatorship was, by far, the longest of political stability (if such a word can be used in the context of Soviet politics).  But on the other hand, the aforementioned stability was also often a period of ideological glaciation so that the period preceding can provide a greater variety of information on the tendencies which reigned among intellectuals trying to adapt their work to the new circumstances. This is what we tried to perceive while studying articles dealing with probability among the first volumes of the encyclopedia.

**Bibliography**

(Borel, 1906)    Borel, Emile. *Avant-Propos*. **Revue du Mois**, 1, 1906

(Bailes, 1978)   Bailes, Kendall E. *Technology and society under Lenin and Stalin : origins of the Soviet technical intelligentsia, 1917-1941*. Princeton University Press, Princeton, 1978

(Bennett-Hurst and Bowley, 1915)    Bennett-Hurst, A.R. and Bowley, Arthur L.. *Livelihood and Poverty: a study in the economic conditions of working-class households*. G. Bell, London, 1915

(Bernstein, 1927)    Бернштейн С.Н. *Теория вероятностей*. Гос. издательство, 1927

(Betti, 2009)    Betti, Renato. *La matematica come abitudine del pensiero - Le idee scientifiche di Pavel Florenskij*. Università Bocconi Centro PRISTEM, 2009

(Borel, 1923)    Борель Э. *Случай. Перевод с французского Ю.И. Костицыной под редакцией В.А.Костицына*. Современные проблемы естествознания, 8. Госиздат, Москва.-Петроград, 1923.

(Bowley, 1907)    Bowley, Arthur. *Elements of Statistics*. London, P. S. King & son; 1907

(Bugaiev, 1898)    Bugaiev, N. V. *Les mathématiques et la conception du monde du point de vue scientifique*. In Rudio, F., editor, Verhandlungen des ersten internationalen Mathematiker-Congresses in Zürich. Teubner, Leipzig., 1898. p. 206–223

(Bustamante, Cléry and Mazliak, 2015)    Bustamante M-C., Cléry M. et Mazliak L. *Le Traité du calcul des probabilités et de ses applications: étendue et limites d'un projet borélien de grande envergure (1921-1939)*. **North-Western European Journal of Mathematics**, 1, 2015. p.85-123

(Castelnuovo, 1919)    Castelnuovo, Guido. *Calcolo delle probabilità*. D. Alighieri, Milano, 1919.

(Czuber, 1914)    Emanuel Czuber. *Wahrscheinlichkeitsrechnung und ihre Anwendung auf Fehlerausgleichung, Statistik und Lebensversicherung*. Teubner, Leipzig, 1903. Second edition 1910, third 1914.

(David-Fox, 1997)    David-Fox, Michael. *Revolution of he Mind: Higher learning among the Bolsheviks*, 1918-1929. Cornell University Press, Ithaca, 1997

(David-Fox, 1998)    Davif-Fox, Michael. *Symbiosis to Synthesis: the Communist Academy and the bolshevization of the Russian Academy of Sciences, 1918-1929*. **Jarbücher für Geschichte Osteureopas**, 46, 1998. p.219-243 (Dubnov and Rashevskij, 1949)   Дубнов Я.С., Рашевский П.К. *В.Ф. Каган*. **Труды семинара по векторному и тензорному анализу**. 7, 1949. p.16-30

(Durand and Mazliak, 2011)   Durand, A. and Mazliak, L. (2011). *Revisiting the sources of Borel's interest in probability : Continued fractions, social involvement, Volterra's prolusione*. **Centaurus**, 53 :306–332.

(Engels, 1959)   Engels, Friedrich *Anti-Dühring (1878)*. Foreign Languages Publishing House, Moscow, 1959

(Ekelund and Hebert, 1999)  Ekelund, Robert B. Jr., Hebert, Robert F.. *Secret Origins of Modern Microeconomics: Dupuit and the Engineers*, University of Chicago Press, 1999.

(Ford, 1991)   Ford, Charles E., 1991. *Dmitrii Egorov: Mathematics and religion in Moscow*. Mathematical Intelligencer 13, 24–30.

(Gliko, 2011)   Глико, А.О. *Отто Юльевич Шмидт в истории России XX века и развитие его научных идей*. Физматлит, Москва, 2011

(Graham, 1987)  Graham, Loren R.   *Science, Philosophy, and Human Behavior in the Soviet Union* Columbia Univ Pr; 1987

(Guieu, 1998)  Guieu, Jean-Michel. *L'engagement européen d'un grand mathématicien français : Émile Borel et la "coopération européenne", des années vingt aux années quarante,* **Bulletin de l'Institut Pierre Renouvin**, n°5, été 1998.

(Hungerford, 2008)  Hungerford, Thomas W. *Algebra*. Graduate Texts in Mathematics; Bd. 73. Springer, New York 2008

(Hayek, 1944)  Hayek, F.A.  *The road to serfdom.*  Routledge and University of Chicago Press, 1944.

(Jaurès, 1901)   Jaurès, Jean. *La philosophie de Vaillant*. **La Petite république socialiste**, 8 janvier 1901.

(Krementsov, 1997)  Krementsov, Nikolai. *Stalinist Science*. Princeton: Princeton University Press, 1997

(Kazanin, 2007)   Казанин Игорь Евгеньевич. *Формирование руководством РСФСР-СССР партийно-государственной политики по отношению к интеллигенции в Октябре 1917-1925 г.* **Диссертация (**Отечественная история**)**. Волгоград 2007.

(Kolmogorov, 1931)  Kolmogoroff Andrei. *Über die analytischen Methoden in der Wahrscheinlichkeitsrechnung*. **Mathematische Annalen** 104, 149-160, 1931.

(Kolmogorov, 1933)   Kolmogoroff Andrei.*Grundbegriffe der Wahrscheinlichkeitsrechnung*, Ergebnisse der Mathematik und ihrer Grenzgebiete, II 3,   Springer, Berlin 1933.

(Koriakin, 2011)   Корякин, Владислав. *Отто Шмидт*. Великие Исторические Персоны. Вече, Москва, 2011.

(Kuznetsov, 2006)	Кузнецов, И.В. *История отечественной журналистики (1917-2000)* Москва, Изд.а Флинта Наука, 2006.

(Khinchin, 1926)	Хинчин, А.Я. *Идеи интуиционизма и борьба за предмет в современной математике*. **Вестник Коммунистической Академии** 16:184–192, 1926.

(Lakhtin, 1924) Лахтин, Л. К. *Курс теории вероятностей. Специальные пособия для высшей школы*. Государственное издательство, Москва-Петроград, 1924.

(Lenin, 1961)	Lenin, Vladimir Ilich *Philosophical Notebooks*. Collected Works, vol. XXXVIII. Foreign Languages Publishing House, Moscow, 1961.

(Lenin, 1918)	Lenin Vladimir Ilich *The Immediate Task of the Soviet Government*. Article on April 28, 1918 in **Pravda** No. 83. In Lenin's Collected Works, 4th English Edition, Progress Publishers, Moscow, 1972 Volume 27, p.235-77

(Lévy, 1925) Lévy, Paul. *Calcul des Probabilités*. Gauthier-Villars, 1925

(Lopshitz and Rashevskij, 1969)	Лопшиц А.М. Рашевский П.К. *Вениамин Федорович Каган (1869–1953)*. Замечательные ученые Московского университета, 39. Издательство Московского Университета, 1969.

(Maksimov, 1931) Максимов, А.А. **За марксистско ленинское естествознание**, 1, 1931.

(Markov, 1924) Марков, Андрей А. Исчисление вероятностей. - 4-е посмертное изд., перераб. автором / биогр. очерк А.С. Безиковича. Госиздат, Москва, 1924.

(Mazliak and Perfettini, 2016)	Mazliak, Laurent et Perfettini, Thomas. *Réflexions sur le hasard dans la société soviétique en construction. Deux textes tirés du débat des années 1920*.To appear in "**Revue de Synthèse**", 2016.

(Mazliak and Sage, 2014)	Mazliak, Laurent and Sage, Marc. *Au-delà des réels : Émile Borel et l'approche probabiliste de la réalité*. **Revue d'histoire des sciences**, 67(2), 2014. p.331–357

(Malia, 1999)	Malia, Martin. *La tragédie soviétique*. Points Histoire. Seuil, Paris. 1999

(Mespoulet, 2001)	Mespoulet, Martine. *Statistique et révolution en Russie. Un compromis impossible (1880-1930)*. Rennes, Presses Universitaires de Rennes, 2001

(Ossorgin, 1955)	Осоргин М. А. *Времена*. Париж, 1955

(Orlov, 1925)	Орлов, Иван Е. *Логика естествознания*. Государственное издательство, Москва, 1925

(Petrov, 1960) Петров, Ф.Н. *Первые советские энциклопедии*. **Книга исследования и материалы**. Сб.3. Издательство всесоюзной книжной палаты, Москва, 1960. p.132-138 (Poincaré, 1912) Poincaré, Henri. *Calcul des Probabilités (2ème édition)*. Gauthier-Villars, 1912.

(Ruud, 1990) Ruud, Charles. *Russian Entrepreneur: Publisher Ivan Sytin of Moscow, 1851-1934.* McGill-Queen's University Press, 1990.

(Rikun, 2012) Рикун, И. Э. *«МАТЕЗИС»—лучшее российское научно-просветительское издательство первой четверти xx века.* **Математика в высшем образовании** 10, 2012. p.141-147

(Selsam and Martel, 1963) Selsam, Howard, and Martel, Harry. *Reader in Marxist philosophy: from the writings of Marx, Engels, and Lenin.* New York: International Publishers, 1963.

(Smit-Falker, 1927) Смит-Фалкер М. *Учёт, статистика и плановость*. **Вестник статистики**, 4, 1927, 15-30.

(Schmidt, 2015) Шмидт, Отто Ю. *Абстрактная теория групп*. Физико-математическое наследие: математика (история математики). Ленанд, Москва, 2015.

(Schmidt, 1922) Шмидт, Отто Ю. *Госиздат за 5 лет*. Госиздат, Москва, 1922.

(Seneta, 2003) Seneta, Eugene. Statistical Regularity and Free Will: L.A.J. Quetelet and P.A. Nekrasov. **International Statistical Review**. 71, 2, 2003. 319-334

(Seneta, 2004) Seneta, Eugene. *Mathematics, religion, and Marxism in the Soviet Union in the 1930s*. **Historia Mathematica** 31, 2004. p.337–367

(Sharansky, 1988) Sharansky, N. *Fear no evil: The classic memoir of one man's triumph over a police state*. New York: Random House, 1988.

(Starzyk, Blatz and Ross, 2009) Starzyk Katherine B., Blatz Craig W., and Ross, Michael. *Acknowledging and Redressing Historical Injustices*. in John T. Jost, Aaron C. Kay, and Hulda Thorisdottir (eds). Social and Psychological Bases of Ideology and System Justification. Oxford University Press, 2009. p.463-479

(Tagliagambe, 2003) Tagliagambe Silvano. *Mathematics and Culture in Russia*, in M. Emmer (editor), *Mathematcs and Culture*, Sprinter-Verlag, New-York, 2003. p. 67-94

(Verburgt, 2016) Verburgt, Lukas *On A. Ia. Khinchin's "Ideas of intuitionism and the struggle for a subject matter in contemporary mathematics"* (1926). To appear in **Historia Mathematica**, 2016.

(Vucinich, 2000) Vucinich, A. *Soviet mathematics and dialectics in the Stalin era*. **Historia Mathematica** 27, 2000. 54–76.

(Werth, 2012) Nicolas Werth. *Histoire de l'Union Soviétique*. Quadrige Manuels. Presses Universitaires de France. Paris, 2012.

(Yanovskaya, 1931) Яновская, София А. **За марксистско ленинское естествознание**,